\documentclass [12pt]{article}
\usepackage{mathrsfs}
\usepackage{amsthm,amsmath,amsfonts,amssymb,graphicx,graphics,psfrag}

\newtheorem{thm}{Theorem}

\newtheorem{lem}{Lemma}

\date{}

\begin{document}

\title{An Improvement on Vizing's Conjecture}

\author{ Yunjian Wu\thanks{Corresponding author: y.wu@seu.edu.cn}
\\ {\small Department of Mathematics}
\\ {\small Southeast University, Nanjing, 211189, China}
}

\maketitle

\begin{abstract}
Let $\gamma(G)$ denote the domination number of a graph $G$. A {\it
Roman domination function} of a graph $G$ is a function $f:
V\rightarrow\{0,1,2\}$ such that every vertex with 0 has a neighbor
with 2. The {\it Roman domination number} $\gamma_R(G)$ is the
minimum of $f(V(G))=\Sigma_{v\in V}f(v)$ over all such functions.
Let $G\square H$ denote the Cartesian product of graphs $G$ and $H$.
We prove that $\gamma(G)\gamma(H) \leq \gamma_R(G\square H)$ for all
simple graphs $G$ and $H$, which is an improvement of
$\gamma(G)\gamma(H) \leq 2\gamma(G\square H)$ given by Clark and
Suen \cite{CS}, since $\gamma(G\square H)\leq \gamma_R(G\square
H)\leq 2\gamma(G\square H)$.
\begin{flushleft}
\end{flushleft}
\end{abstract}

\vskip 1mm \vskip 1mm

\section{Introduction}

In this note, we consider simple finite graphs only and follow
\cite{hay} for terminology and definitions.

let $G = (V, E)$ be a graph with vertex set $V$ and edge set $E$.
For any vertex $v\in V$, the {\it open neighborhood of $v$} is the
set $N(v)=\{u\in V\ |\ uv\in E\}$ and the {\it closed neighborhood}
is the set $N[v]=N(v)\cup\{v\}$. For a set $S\subseteq V$, the open
neighborhood is $N(S)=\bigcup_{v\in S}N(v)$ and the closed
neighborhood is $N[S]=N(S)\cup S$. A set $S\subseteq V$ is a {\it
dominating set} of $G$ if every vertex not in $S$ is adjacent to a
vertex in $S$. The {\it domination number} of $G$, denoted by
$\gamma(G)$, is the minimum cardinality of a dominating set. A
domination set of cardinality $\gamma(G)$ is called a {\it
$\gamma$-set} of $G$. Recently, a variant of the domination
number---Roman domination number is suggested by Stewart \cite{Ste}.
A {\it Roman dominating function} (RDF) on a graph $G=(V,E)$ is a
function $f: V\rightarrow\{0,1,2\}$ satisfying the condition that
every vertex $u$ for which $f(u)=0$ is adjacent to at least one
vertex $v$ for which $f(u)=2$. The weight of $f$ is
$f(V(G))=\Sigma_{v\in V}f(v)$. The {\it Roman domination number},
denoted by $\gamma_R(G)$, equals the minimum weight of an RDF of
$G$, and we say that a function $f$ is a {\it
$\gamma_R(G)$-function} if it is an RDF and $f(V(G))=\gamma_R(G)$.
For a graph $G$, let $f: V\rightarrow\{0,1,2\}$, and let $(V_0, V_1,
V_2)$ be the order partition of $V$ induced by $f$, where
$V_i=\{v\in V(G)\ |\ f(v)=i\}$ for $i=0,1,2$. Note that there exists
a 1-1 correspondence between the functions $f:
V\rightarrow\{0,1,2\}$ and the ordered partitions $(V_0, V_1, V_2)$
of $V(G)$. Thus we will write $f=(V_0, V_1, V_2)$.

\vspace{4mm}

Cockayne et al. \cite{Coc} showed the following results.

\begin{lem}\label{inequ}
{\rm{(\cite{Coc})}} For any graph $G$, $\gamma(G)\leq
\gamma_R(G)\leq 2\gamma(G)$.

\end{lem}

\begin{lem}\label{roman}
{\rm{(\cite{Coc})}} Let $f=(V_0, V_1, V_2)$ be any
$\gamma_R(G)$-function. Then $V_2$ is a $\gamma$-set of $G[V_0\cup
V_2]$.
\end{lem}

For a pair of graphs $G$ and $H$, the Cartesian product $G\square H$
of $G$ and $H$ is the graph with vertex set $V(G)\times V(H)$ and
where two vertices are adjacent if and only if they are equal in one
coordinate and adjacent in the other. In 1963, V.\ G.\ Vizing
\cite{VGV} conjectured the following:

\vspace{2mm}

\noindent{\bf Vizing's Conjecture}. For any graphs $G$ and $H$,
$\gamma(G)\gamma(H) \leq \gamma(G\square H)$.

\vspace{2mm}

We note that there are graphs $G$ and $H$ for which the above
equality holds. The reader is referred to Hartnell and Rall
\cite{HR} for a  summary of recent progress on Vizing's conjecture.
Recently, Clark and Suen \cite{CS} gave the following result.

\begin{thm}\label{2times}
{\rm{(\cite{CS})}} For any graphs $G$ and $H$, $\gamma(G)\gamma(H)
\leq 2\gamma(G\square H)$.
\end{thm}

We shall show in this note that $\gamma(G)\gamma(H) \leq
\gamma_R(G\square H)$, which is an improvement of
$\gamma(G)\gamma(H) \leq 2\gamma(G\square H)$ by Lemma \ref{inequ}.

\section{Main results}

\begin{thm}
For any graphs $G$ and $H$,
$$\gamma(G)\gamma(H) \leq \gamma_R(G\square H).$$
\end{thm}

\noindent {\bf Proof.}  Let $f=(V_0, V_1, V_2)$ be any
$\gamma_R(G\square H)$-function of graph $G\square H$. Denote
$D=V_1\cup V_2$. By Lemma \ref{roman}, $D$ and $V_2$ are domination
set of graphs $G\square H$ and $G\square H-V_1$, respectively. Let
$\{u_1, u_2, \ldots, u_{\gamma(G)}\}$ be a dominating set of $G$.
Then we partition $V(G)$ into $\gamma(G)$ sets $\{\Pi_1, \Pi_2,
\ldots, \Pi_{\gamma(G)}\}$ satisfying the following properties:\\
\vspace{-2mm}

$($i$)$ $u_i \in \Pi_i$, \\

\vspace{-2mm}

$($ii$)$ $u \in \Pi_i$ implies $u= u_i$ or $u$ is adjacent to~$u_i$.

Note that this partition is not unique. The partition of $V(G)$
induces a partition $\{D_1, D_2, \ldots, D_{\gamma(G)}\}$ of $D$
where
$$
D_i = (\Pi_i \times V(H)) \cap D.
$$
Let $P_i$ be the projection of $D_i$ onto $H$.  Then
$$
P_i = \{v \,\vert\, (u,v) \in D_i \mbox{\rm{ for some }} u \in
\Pi_i\}.
$$
For any $i$, $P_i \cup (V(H) - N_H[P_i])$ is a dominating set of
$H$, so the number of vertices in $V(H)$ not dominated by $P_i$
satisfies the inequality
\begin{eqnarray} \label{lowerbound}
|V(H) - N_H[P_i]| \geq \gamma(H) - |P_i|.
\end{eqnarray}
For $v \in V(H)$,  denote
$$
Q_v = V_2 \cap (V(G) \times \{v\}) = \{ (u,v) \in V_2 \,\vert\, u
\in V(G)\},
$$
let $C$ be the subset of $\{1,2,\ldots,\gamma(G)\} \times V(H)$
given by
$$
C = \{\,  (i,v) \, \vert \,  \Pi_i \times \{v\} \subseteq
N_{G\square H}[ \, Q_v \, ] \, \}.
$$
Set
\begin{eqnarray*}
L_i & = & \{ (i,v) \in C \, \vert \, v\in V(H) \}, \\
R_v & = &\{ (i,v) \in C \, \vert\,  1 \le i \le \gamma(G) \}.
\end{eqnarray*}
It is clear that
$$
N = \vert C \vert = \sum _{i=1}^{\gamma(G)}\vert L_i\vert =  \sum
_{v \in V(H)} \vert R_v \vert .
$$
If $v \in V(H) - N_H[P_i]$, then the vertices in $\Pi_i \times
\{v\}$ must be dominated by vertices in $Q_v$ since $\Pi_i \times
\{v\}\nsubseteq D$ and $V_2$ is a dominating set of graph $G\square
H-V_1$. Therefore $(i,v) \in L_i$. This implies that $|L_i| \ge
|V(H) - N_H[P_i]|$. Hence
$$
N \geq \sum_{i=1} ^{\gamma(G)} |V(H) - N_H[P_i]|
$$
Now it follows from (\ref{lowerbound})  that
\begin{eqnarray*}
N &\geq& \gamma(G)\gamma(H) - \sum_{i=1}^{\gamma(G)} |P_i|  \\
&\geq& \gamma(G)\gamma(H) - \sum_{i=1}^{\gamma(G)} |D_i|.
\end{eqnarray*}
So we obtain the following lower bound for $N$.
\begin{eqnarray}  \label{first}
N \geq \gamma(G)\gamma(H) - |D|=\gamma(G)\gamma(H)-|V_1|-|V_2|.
\end{eqnarray}
For each $v \in V(H)$, $|R_v| \leq |Q_v|$. If it is not true, then
$$
\{\, u \,\vert\, (u,v) \in Q_v\} \cup \{ \, u_j \,  | \,  (j,v)
\notin R_v \}
$$
is a dominating set of $G$ with cardinality
$$
|Q_v| + (\gamma(G) - |R_v|) = \gamma(G) - (|R_v|-|Q_v|) < \gamma(G),
$$
and we have a contradiction.  This observation shows a upper bound
for $N$.
\begin{equation}
\label{second} N = \sum_{v \in V(H)}  |R_v| \le \sum_{v \in V(H)}
|Q_v| = \vert V_2 \vert.
\end{equation}
It follows from (\ref{first}) and (\ref{second}) that
$$
\gamma(G)\gamma(H)-|V_1|-|V_2| \leq N \le |V_2|,
$$

So we get $\gamma(G)\gamma(H) \leq |V_1|+2|V_2|=\gamma_R(G\square
H)$.\hfill$\Box$

\vspace{2mm}

\noindent{\bf Remark:} One may wonder if there is a similar result
on Roman domination number as Vizing's conjecture. In fact, there
are examples showing the inequality $\gamma_R(G)\gamma_R(H) \leq
\gamma_R(G\square H)$ fails, e.g., $\gamma_R(K_2)=2$, but
$\gamma_R(K_2\square K_2)=3$.

\vskip 20pt

\noindent \title{\Large\bf Acknowledgments} \maketitle
 The authors are indebted to Professor Qinglin Yu for his reading of the manuscript and the constructive comments.

\end{document}